\documentclass{article}

\usepackage{amsmath,amsfonts,amssymb}
\usepackage[colorlinks=true,citecolor=blue]{hyperref}	
\usepackage[normalem]{ulem}
\renewcommand{\o}{\omega}

\newcommand{\bC}{\mathbb{C}}  
\newcommand{\bR}{\mathbb{R}}  
  
\newcommand{\bH}{\mathbb{H}}

\newcommand{\ga}{\mathfrak{a}}  
\newcommand{\gb}{\mathfrak{b}}  
\newcommand{\gc}{\mathfrak{c}}

\renewcommand{\gg}{\mathfrak{g}}  
\newcommand{\gh}{\mathfrak{h}}  
  
\newcommand{\gl}{\mathfrak{l}}  
\newcommand{\gm}{\mathfrak{m}}

\newcommand{\gs}{\mathfrak{s}}  
\newcommand{\gt}{\mathfrak{t}}  
\newcommand{\gu}{\mathfrak{u}}

\newcommand{\so}{\mathfrak{so}}  
\newcommand{\su}{\mathfrak{su}}  
  
\newcommand{\gsp}{\mathfrak{sp}}  
\newcommand{\gsl}{\mathfrak{sl}}  
\newcommand{\ggl}{\mathfrak{gl}}

\newcommand{\ra}{\rightarrow}

\newtheorem{theorem}{Theorem}
\newtheorem{Th}{Theorem} 
\newtheorem{Prop}{Proposition}  
\newtheorem{Cor}{Corollary}  
\newtheorem{Lem}{Lemma}  
\newtheorem{Def}{Definition}  
  
\newtheorem{Rm}{Remark}

\newcommand{\bt}{\begin{Th}}  
\newcommand{\et}{\end{Th}}  
\newcommand{\bp}{\begin{Prop}}  
\newcommand{\ep}{\end{Prop}}  
\newcommand{\bc}{\begin{Cor}}  
\newcommand{\ec}{\end{Cor}}  
\newcommand{\bl}{\begin{Lem}}  
\newcommand{\el}{\end{Lem}}  
\newcommand{\bd}{\begin{Def}}  
\newcommand{\ed}{\end{Def}}  
\newcommand{\bpf}{\begin{proof}}  
\newcommand{\epf}{\end{proof}} 

\newcommand{\n}{\nabla}  
  
\newcommand{\ot}{\otimes}

\newcommand{\be}{\begin{equation}}  
\newcommand{\ee}{\end{equation}}  
  
\newcommand\re[1]{(\ref{#1})}  
\newcommand{\arr}{\begin{array}{rlll}}  
\newcommand{\ea}{\end{array}}  
\newcommand{\bea}{\begin{eqnarray}}  
\newcommand{\eea}{\end{eqnarray}}  
\newcommand{\bean}{\begin{eqnarray*}}  
\newcommand{\eean}{\end{eqnarray*}}  

\newcommand{\im}{\sqrt{-1}\,}

\def\address#1#2{\begingroup
\noindent\parbox[t]{7.8cm}{
\small{\scshape\ignorespaces#1}\par\vskip1ex
\noindent\small{\itshape E-mail address}
\/: #2\par\vskip4ex}\hfill
\endgroup}

\title%[Unimodular Sasaki and Vaisman Lie groups]
{Unimodular Sasaki and Vaisman Lie groups
}

\author{Vicente Cort\'es 
and Keizo Hasegawa
}

\begin{document}  

\maketitle

\let\thefootnote\relax\footnotetext{
 \par
The research of V.C.\ was partially funded by the Deutsche Forschungsgemeinschaft (DFG, German Research Foundation) under Germany's Excellence Strategy -- EXC 2121 Quantum Universe  -- 390833306.
\par
The research of K.H.\ was supported by JSPS (Japan Society for the Promotion of Science) KAKENHI Grant Number 21K03248.
} 

\begin{abstract}  
This is a continuation of our study on
homogeneous locally conformally K\"ahler and Sasaki manifolds.
In a recent work \cite{AHK}, applying the technique of modification we have determined
all homogeneous Sasaki and Vaisman manifolds of unimodular 
Lie groups, up to modifications. In this paper we determine all such modifications explicitly 
for the case of Lie groups, obtaining a complete 
classification of unimodular Sasaki and Vaisman Lie groups. Furthermore,
we determine the biholomorphism type of a simply connected unimodular 
Vaisman Lie group of each type. 
\end{abstract}  

%%%%%%%%%%%%%%%%%%%%%%%%%%%%%%  NEWSECTION %%%%%%%%%%%%%%%%%%%%%%
%%%%%%%%%%%%%%%%%%%%%%%%%%%%%%  NEWSECTION %%%%%%%%%%%%%%%%%%%%%%
\section*{Introduction}

In a series of our papers (\cite{ACHK}, \cite{AHK}, \cite{HK1}, \cite{HK2}) 
we have studied homogeneous locally conformally K\"ahler (shortly lcK)
and Sasaki manifolds. In our recent paper \cite{AHK}, applying the technique of modification we have determined
all homogeneous Sasaki and Vaisman manifolds of unimodular Lie groups, up to modifications. In particular,
for the case of unimodular Sasaki and Vaisman Lie groups we have obtained the following basic result.

%{\bf Theorem} \quad
\begin{theorem}
{\em
A Sasaki unimodular Lie algebra is, up to modification, isomorphic to one of the three types:
$\gh_{2m+1}$, $\su(2)$, $\gsl(2, \bR)$, where $\gh_{2m+1}$ is the Heisenberg Lie algebra of dimension $2m+1$.
Accordingly, a Vaisman unimodular Lie algebra is, up to modification,
isomorphic to one of the following:
$$\bR \oplus \gh_{2m+1},\, \bR \oplus \su(2),\,\bR \oplus \gsl(2, \bR). \eqno{(*)}$$
In terms of Lie groups, a simply connected Sasaki unimodular Lie group is, up to modification,
isomorphic to one of the three types: $H_{2m+1}, SU(2), \widetilde{SL}(2, \bR)$,
where $H_{2m+1}$ is the Heisenberg group of dimension $2m+1$.
Accordingly, a simply connected Vaisman unimodular Lie group is, up to modification,
isomorphic to one of the following:
$$\bR \times H_{2m+1},\, \bR \times SU(2),\, \bR \times \widetilde{SL}(2, \bR). \eqno{(**)}$$
}
\end{theorem}

Note that we have $\gu(2) = \bR \oplus \su(2),\,\gg\gl(2, \bR)= \bR \oplus \gsl(2, \bR)$; and from now on
we denote the nilpotent Lie algebra
$\bR \oplus \gh_{2m+1}$  by \,$\gg\gh_{2m+2}$ or $\gg\gh$ for short, and the corresponding Lie group $\bR \times H_{2m+1}$ by
$GH_{2m+2}$ or $GH$ for short.

In this paper we determine all possible modifications of Sasaki and Vaisman Lie algebras in the list of
the above theorem; and thus obtain a complete classification of simply connected Sasaki and Vasiman
Lie groups. It turns out that any modification of a reductive Hermitian Lie algebra with $1$-dimensional center
(in particular, $\bR \oplus \su(2),\, \bR \oplus \gsl(2, \bR)$) is isomorphic to the original one, whereas any modification $\gg\gh_{2m+2}'$ of
$\gg\gh_{2m+2}$ preserves the Vaisman structures, defining a new family of unimodular Vaisman solvable 
Lie algebras (see Section \ref{unimodular}). Similarly, we can see that any modification of a semi-simple Sasaki Lie algebra with trivial center is
isomorphic to the original one, whereas any modification $\gh_{2m+1}' $ of $\gh_{2m+1}$ preserves the Sasaki structure,
defining a new family of unimodular Sasaki solvable Lie algebras (see Section \ref{unimodular}). 
Note that some partial results on unimodular Sasaki and/or Vaisman Lie groups (from different perspectives) have been known lately, 
cf. \cite{AFV}, \cite{AO}, \cite{BK}, \cite{SW2}; and for non-unimodular Sasaki and/or Vaisman Lie groups, e.g. \cite{AFV}, \cite{LL}.

\vspace{10pt}

In Section \ref{complex_str} we determine all the complex structures on $\gu(2), \gg\gl(2, \bR)$ and $\gg\gh_{2m+2}$,
partially extending and elaborating the results
in our previous papers. 
We give unified proofs for the cases of $\gu(2), \gg\gl(2, \bR)$ based on the fact that their complexifications are the same.  
It appears that the determination of all complex structures on $\gg\gh_{2m+2}$ is not known; 
although it can be obtained in a similar way as for the cases of $\gu(2), \gg\gl(2, \bR)$.
In Section \ref{Vaisman_str} we determine all Vaisman structures compatible with the complex structures on $\gu(2), \gg\gl(2, \bR)$ and $\gg\gh_{2m+2}$
which are determined in Section \ref{complex_str}. We again make the proofs for the cases of $\gu(2), \gg\gl(2, \bR)$ in a unified form,
but the consequences are quite different: as we have seen in our previous papers, while all the lcK
structures on $\gu(2)$ are of Vaisman type, there are on $\gg\gl(2, \bR)$
both lcK structures of Vaisman type and those of  non-Vaisman type.
For the case of $\gg\gh_{2m+2}$, we know (\cite{SW} \cite{HK1}) that the only lcK nilpotent Lie algebras
are of Heisenberg type, that is, $\gg\gh_{2m+2}$; and all these lcK structures are of Vaisman type. Here, 
we determine all Vaisman structures on  $\gg\gh_{2m+2}$ compatible with the complex structures we have obtained in Section \ref{complex_str}.
In Section \ref{modification}, we review and discuss in some detail the technique of modification in the category of Hermitian 
Lie algebras and Sasaki Lie algebras. The key points are that Modification is equivalence relation and
preserves unimodularity, as shown in our previous paper. In section \ref{unimodular} we determine all possible modifications of 
$\gu(2), \gg\gl(2, \bR)$ and $\gg\gh_{2m+2}$; and thus obtain a complete classification of unimodular Sasaki and Vaisman 
Lie algebras. In Section \ref{biholomorphism}, we see that the complex structures on $\gu(2), \gg\gl(2, \bR)$ and $\gg\gh_{2m+2}$ 
are biholomorphic to $\bC^2 \backslash \{0\}, \bC \times \bH$ and $\bC^{m+1}$ respectively; and 
any modification $\gg\gh_{2m+2}'$ of $\gg\gh_{2m+2}$ preserves the original complex structure $\bC^{m+1}$.
%%%%%%%%%%%%%%%%%%%%%%%%%%% NEWSECTION %%%%%%%%%%%%%%%%%%%%%%
%%%%%%%%%%%%%%%%%%%%%%%%%%% NEWSECTION %%%%%%%%%%%%%%%%%%%%%%

\section{Complex structures on $\mathfrak{gl}(2,\bR), \gu(2)$ and $\gg\gh_{2m+2}$} \label{complex_str}
Recall that any left-invariant complex structure on a Lie group $G$ with Lie algebra $\mathfrak g$ induces an endomorphism 
$J\in \mathrm{End}\, \mathfrak{g}$ such that $J^2=-\mathrm{Id}$ and $N_J=0$, where $N_J\in \bigwedge^2 \mathfrak{g}^*\otimes \mathfrak{g}$
is the Nijenhuis tensor of $J$, defined with the help of the Lie bracket of $\mathfrak g$ by
\be \label{Nj:eq} N_J (X,Y) := [JX,JY] -[X,Y] -J[X,JY]-J[JX,Y] \ee
for $X,Y\in \gg$.

An endomorphism $J\in \mathrm{End}\, \mathfrak{g}$ such that $J^2=-\mathrm{Id}$ and $N_J=0$ is called a {\em complex structure} on the Lie algebra $\mathfrak g$ and defines a left-invariant
complex structure on any Lie group with Lie algebra $\mathfrak g$. 

\bd
Complex structures $J_1$ and $J_2$ on a Lie algebra $\gg$ are equivalent if there exists an automorphism $\phi$ of
$\gg$ such that $\phi J_1 = J_2 \phi$.
\ed

In this section we determine all complex structures on $\gu(2), \gg\gl(2, \bR)$ and $\gg\gh_{2m+2}$,
partially extending and elaborating the arguments in our previous papers. 

We already know all complex structures on 
$\gg\gl(2,\bR)=\bR \oplus \mathfrak{sl}(2, \bR)$ and $\gu(2)=\bR \oplus \mathfrak{su}(2)$ (\cite{ACHK}, \cite{HK1}, \cite{SS}).
Since the complexifications of both Lie algebras are the same complex Lie algebra $\gg\gl(2, \bC)$, 
we can treat the determination of all complex structures on them simultaneously.

For  $\gg_1 := \bR \oplus \mathfrak{sl}(2, \bR)$, 
take a basis $\{X_1, Y_1, Z_1 \}$ for $\mathfrak{sl}(2, \bR)$
with Lie bracket defined by
\be \label{sl2:eq} [X_1,Y_1]=Z_1,\;[Z_1,X_1]=-Y_1,\, [Z_1, Y_1]=X_1,\ee
and $T$ as a generator of the center $\bR$ of $\gg_1$, where we set
$$ 
X_1=\frac{1}{2}\left(
\begin{array}{cc}
0 & 1\\
1 & 0
\end{array}
\right),
\;
Y_1=\frac{1}{2}\left(
\begin{array}{cc}
1 & 0\\
0 & -1
\end{array}
\right),
\;
Z_1=\frac{1}{2}\left(
\begin{array}{cc}
0 & -1\\
1 & 0
\end{array}
\right).
$$

For  $\gg_2 := \bR \oplus \mathfrak{su}(2)$, take a basis $\{X_2, Y_2, Z_2 \}$
with Lie bracket defined by
\be \label{su2:eq} [X_2,Y_2]=-Z_2,\, [Z_2, X_2]=-Y_2,\, [Z_2,Y_2]=X_2, \ee
and $T$ a generator of the center of $\gg_2$, where we set
$$ 
X_2=\frac{1}{2}\left(
\begin{array}{cc}
0 & -1\\
1 & 0
\end{array}
\right),
\;
Y_2=\frac{\sqrt{-1}}{2}\left(
\begin{array}{cc}
0 & 1\\
1 & 0
\end{array}
\right),
\;
Z_2=\frac{\sqrt{-1}}{2}\left(
\begin{array}{cc}
1 & 0\\
0 & -1
\end{array}
\right).
$$

\vspace{10pt}

We know that there exists a one-to-one correspondence between complex structures $J$ on a (real) Lie algebra
$\gg$ and complex Lie subalgebras $\gh \subset \gg_{\bC}$ which satisfy $\gg_{\bC} = \gh \oplus \overline{\gh}$
(direct sum of vector spaces). It is given by
$J \mapsto \gh_J :=\{X- \sqrt{-1} JX | X \in \gg\}$. 

For $\gg_1$ and $\gg_2$,
we have the same complexification $\gg_{\bC}$ of $\gg_i \; (i=1,2)$ given by
$$\gg_{\bC}= \mathfrak{gl}(2, \bC) = \mathfrak{gl}(2, \bR) \otimes \bC= \mathfrak{u}(2) \otimes \bC.$$
We determine all complex subalgebras $\gh$ of $\gg_{\bC}$ satisfying
$\gg_{\bC} = \gh \oplus \overline{\gh}$, and thus determine all integrable
complex structures on $\gg_i \;(i=1,2)$. Note that the complex conjugation is defined
by its fixed point set $\mathfrak{g}_i$.

\vspace{10pt}

We define a basis $\beta_i =\{T, U, V_i, W_i\}, \; i=1,2$, for $\gg_{\bC}$ induced from $\gg_i$,
where $T$ is a generator of the center $\bC$ of $\gg_{\bC}$, and

$$U_i= \im Z_i,\;
V_i=\frac{1}{2}(X_i-\im Y_i),\;
W_i=-\frac{1}{2}(X_i+\im Y_i).$$
The non-zero Lie bracket are given by
$$[V_i, W_i]= (-1)^{i} \frac{1}{2} U_i,\; [U_i, V_i]= V_i,\; [U_i, W_i]=-W_i;$$
and the conjugation with respect to $\gg_i,\; i=1,2$, is given by
$$\overline{T}=T,\; \overline{U_i}=-U_i,\; \overline{V_i}=-W_i,\;
\overline{W_i}=-V_i.$$

In the following, we can safely omit subscripts $i$,
since the Lie brackets are the same except the sign for $[V_i, W_i] =(-1)^{i} \frac{1}{2} U_i$,
and conjugations are the same for $\beta_i,\; i=1,2$.

\vspace{10pt}

Let $\ga$ be the center
of $\gg_{\bC}$ generated by $T$ and $\gb$ the subalgebra of
$\gg_{\bC}$ generated by $U,V,W$; that is, we have
$$\gg_{\bC} = \ga \oplus \gb$$
where $\ga=<T>_{\bC}, \gb=<U, V, W>_{\bC}$.
Let $\pi$ be the projection $\pi: \gg_{\bC} \rightarrow \gb$ and
$\gc$ the image of $\gh$ by $\pi$, then $\gc$  is a subalgebra of $\gb$ which satisfies
$$\gb=\gc+\overline{\gc},$$
and ${\rm dim}\, \gc \cap \overline{\gc} =1$.
We can set a basis $\eta$ of $\gh$ as
$\eta=\{P+Q, R\}\; (P \in \ga \setminus \{ 0\}, Q,R \in \gb)$
%here we use g_C = h + \bar h 
such that $Q \in \gc \cap \overline{\gc}$, and
$\gamma =\{Q,R\}$ is a basis of $\gc$:
$$\gh=<P+Q, R>_{\bC},\; \gc=<Q, R>_{\bC},\; \gb =<Q, R, \overline{R}>_{\bC}.$$
Furthermore, since we have $[Q,R]=[P+Q,R]\in \gc \cap \gh = <R>_{\bC}$ and $[\gb, \gb]=\gb$, we see that
$[Q, R] =\sigma R$ with $\sigma \in \bC\setminus \{ 0\}$; and since $\overline{Q} \in \gc \cap \overline{\gc}$
we can also assume $Q+\overline{Q}=0$
so that $Q$ is of the form $aU+bV+\overline{b}W \, (a \in \bR, b \in \bC)$.

\vspace{10pt}

We first consider the case where $R=qV+rW \, (q,r \in \bC)$. 
We see by simple calculation that if $b \not=0$,
then $q=sb, r=s \overline{b}$
for some non zero constant $s \in \bC$. But then $\displaystyle \overline{R}=-\frac{\bar{s}}{s} R$,
contradicting the fact that $\beta=\{Q, R,\overline{R}\}$ is 
a basis of $\gb$.
Hence we must have $b=0$, and $q \not=0, r =0$ with
$\sigma=a$ or $q =0, r \not=0$ with $\sigma=-a$.
Therefore we can take, as a basis of $\gh$, $\eta=\{T+\delta U, V\}$ or
$\{T+\delta U, W\}$ with $\delta=k+\im l \in \bC$ and $k\neq 0$:
$$\gh=<T+\delta U, V>_{\bC} \;{\rm or}\; <T+\delta U, W>_{\bC}.$$
Note that the latter defines a conjugate complex structure $\overline{\gh}$ of the former type $\gh$
with the parameter $-\bar{\delta}$. 
These conjugate complex structures are not equivalent (unless $l=0$)
but define biholomorphic
complex structures on the associated Lie group $G$ (see Section~\ref{biholomorphism}). 
%(Two complex structures on a Lie algebra are called equivalent if they are related by an automorphism.)

\vspace{10pt}

In the case where $R=pU+qV+rW,\,p, q, r \in \bC$ with $p \not=0$,
we show that there exists an automorphism $\widehat{\phi}$ on $\gg_{\bC}$ 
which maps $\gh_0$ to $\gh$, preserving the conjugation, where
$\gh_{0}$ is a subalgebra of $\gg_{\bC}$ of the first type with $p=0$.
As in the first case, we must have $[Q,R]=\sigma R$, $\sigma \in \bC\setminus \{ 0\}$. We may assume that $p=1$.
We see, by simple calculation that
$$(a-\sigma)q=b,\; (a+\sigma)r=\overline{b},\; (-1)^{i}(b r - \overline{b} q) = 2 \sigma$$
with $b, q, r  \not= 0$ and $a\neq \pm \sigma$, from which we get

$$a^2+(-1)^i |b|^2=\sigma^2,\quad qr=-(-1)^i,$$
$$|q|^2-|r|^2 = \frac{4 a {\rm Re}(\sigma)}{|b|^2}.$$

From the first equation, we see that $\sigma$ must be real or pure imaginary. If $\sigma$ is pure imaginary, the third equation vanishes.
We have then $[R, \overline{R}]=0$, which contradicts the fact that $\gb$ is semi-simple and thus contain no abelian
ideals.

We can define an automorphism $\phi$ on $\gb$ by setting
$$\phi(U)=\frac{1}{\sigma} Q,\, \phi(V)=\frac{|b|}{2 \sigma} R, \, \phi(W)=-\frac{|b|}{2 \sigma} \overline{R}.$$
It extends to an automorphism $\widehat{\phi}$
on $\gg_{\bC}$ which satisfies the required condition. Hence we have the following.

%%%%%%% Proposition %%%%%%%%
\bp 
Each $\gg_i,\, i=1,2$, admits a family of complex structures $J_{i, \delta}, \delta=k~+~\im l$ $(k \not=0)$, 
defined by
$$J_{i,\delta} (T-l Z)=k Z,\; J_{i,\delta} (k Z)=-(T-l Z),\; J_{i,\delta} X_i= \pm Y_i,
\; J_{i,\delta} Y_i= \mp X_i.$$
Conversely, the above family of complex structures exhausts all
complex structures on each of the Lie algebras $\gg_i$.
\ep

%{\color{red}Note that the complex subalgebra $\mathfrak{h}_J\subset \mathfrak{g}_{\bC}$ associated 
%with $J=J_{i,\delta}$ is $< T-\bar{\delta}U,V>_{\bC}$ or $< T-\bar{\delta}U,W>_{\bC}$, 
%respectively (depending on the choice of sign $\pm$).}

We next determine all complex structures on $\gg\gh=\bR \oplus \gh_{2m+1}$ in the same vein as
in the previous cases.
$\gg\gh$ is a nilpotent Lie algebra with a standard basis $\beta =\{T, X_i,Y_j,Z\mid i, j=1,\ldots ,m\}$
for which non-zero Lie brackets are given by
$$[X_i,Y_i]=Z,\; i=1, 2, ..., m.$$
%where $\gh_{2m+1}$ is the Heisenberg Lie algebra.

Let $\gg\gh_{\bC} = \bC \oplus \gh_{{2m+1}, \bC}$ be the complexification of $\gg\gh$,
where $\gh_{{2m+1}, \bC}$ is the complex Heisenberg Lie algebra.
We have a basis $\{T, U, V_i, W_j \mid  i,j =1, \ldots ,m\}$ for $\gg\gh_{\bC}$,
where $T$ is a generator of the center $\bC$ of $\gg\gh_{\bC}$, and
$$U=-\im Z,\;
V_i=\frac{1}{2}(X_i-\im Y_i),\;
W_j=-\frac{1}{2}(X_j+\im Y_j).$$
The non-zero Lie brackets are given by
$$[V_i, W_i]= \frac{1}{2} U,\; i=1,2,...,m;$$
and the conjugation with respect to $\gg$ is given by
$$\overline{T}=T,\; \overline{U}=-U,\; \overline{V_i}=-W_i,\;
\overline{W_j}=-V_j.$$

Let  $\ga=<T>_{\bC}$ be the center
of $\gg\gh_{\bC}$ and $\gb=\gh_{{2m+1}, \bC}=<U, V_i, W_j\mid  i,j =1, \ldots ,m>_{\bC}$ the subalgebra of
$\gg\gh_{\bC}$ generated by $U,V_i,W_j$; that is, we have
$$\gg\gh_{\bC} = \ga \oplus \gb.$$
Let $\pi$ be the projection $\pi: \gg\gh_{\bC} \rightarrow \gb$ and
$\gc$ the image of $\gh$ by $\pi$, then $\gc$  is a subalgebra of $\gb$ which satisfies
$$\gb=\gc+\overline{\gc},$$
and ${\rm dim}\, \gc \cap \overline{\gc} =1$.
We can set a basis $\eta$ of $\gh$ as
$\eta=\{P+Q, R_i\}\; (P \in \ga\setminus \{ 0\}, Q,R_i \in \gb,\; i=1,2,...,m)$
such that $Q \in \gc \cap \overline{\gc}$, and
$\gamma =\{Q,R_i\}$ is a basis of $\gc$:
$$\gh=<P+Q, R_i>_{\bC},\; \gc=<Q, R_i>_{\bC},\; \gb =<Q, R_i, \overline{R}_j>_{\bC}.$$
Furthermore, since $\overline{Q} \in \gc \cap \overline{\gc}$
we can also assume $Q+\overline{Q}=0$
so that $Q$ is of the form $aU+\sum_i b_i V_i+\sum_j \overline{b_j} W_j \,(a \in \bR, b_i \in \bC)$.
%And since ${\rm dim}\, [\gb, \gb] = 1$, we may assume that $[Q, R_i]= \sigma Q$ for some $\sigma \in \bC$, and
%thus we get $b_i = 0, i=1,2,..., m$ and $Q=a U$. Then, we must have $[Q, R_i]=0$ (since $[Q, R_i] \in [\gh,\gh]$ while
%$Q \notin \gh$), $[R_i, \overline{R_i}]= \frac{1}{2} \varepsilon_i Q\;(\varepsilon_i \in \bR, \varepsilon_i \not=0)$ 
%and $[R_i, \overline{R}_j]= 0$ for $i \not= j$.
%We can also assume that $R_k \in \;<V_i, W_j>_{\bC}$ and $\varepsilon_i =\pm 1$.
 Since we have $[\gb, \gb] = <U>_{\bC}$, we can set  $[Q, R_i]= \alpha_i U, [R_j, R_k]= \beta_{j,k} U \;(\alpha_i, \beta_{j,k}  \in \bC)$
 and thus $[Q, \overline{R}_i]= \bar{\alpha_i} U, [\bar{R}_j, \bar{R}_k]= -\bar{\beta}_{j,k} U$. 
It follows that $[\gc, \gc] =\overline{[\gc, \gc]}$, and thus it is a subspace of $\gc \cap \overline{\gc} =<Q>_{\bC}$.
In particular, we get $[Q, R_i]= \sigma_i Q$, $[R_j, R_k]=\tau_{j,k} Q \; (\sigma_i, \tau_{j,k}  \in \bC)$.
Since $[R_j,R_k], [Q, R_i] \in [\gh,\gh]$ while $Q \notin \gh$, we must have $\sigma_i =0, \tau_{j,k} =0$ for all $i,j,k$. 
Hence we have $[Q, R_i]=[R_j,R_k]=0$ for all $i,j,k$; and thus $[\gc, \gc] =\overline{[\gc, \gc]}=0$.
In particular, $Q$ is in the center of %the Heisenberg Lie algebra 
$\mathfrak b$, i.e.\  $Q=a U \,(a \not=0)$ and $b_i=0$.
Therefore, we have $<U>_{\bC} = [\gb,\gb]=[\gc, \overline{\gc}] = <Q>_{\bC}$. 
We can also assume (without loss of generality) the following:
$$a=1,\;[R_i, \overline{R_i}]= \frac{1}{2} \varepsilon_i Q,\;
[R_i, \overline{R}_j]= 0 \;(i \not= j),$$
where $R_k \in \;<V_i, W_j>_{\bC}$ and $\varepsilon_i =\pm 1$.

For $\gb=<Q, R_i, \overline{R}_j>_{\bC}$ we define an automorphism $\phi$ on $\gb$ by
$$\phi(U)=Q, $$
$$\phi(V_i)= \, R_i, \phi(W_i)= -\, \overline{R}_i \quad \mbox{for}\quad \varepsilon_i=-1,$$
$$\phi(V_i)= \, \overline{R}_i, \phi(W_i)= -\, R_i \quad\mbox{for}\quad \varepsilon_i=1,$$
which extends to an automorphism $\widehat{\phi}$ on $\gg\gh_{\bC}$ preserving the conjugation.
We have thus shown the following.

%%%%%%% Proposition %%%%%%%%
\bp \label{heis:prop} The Lie algebra 
$\gg\gh = \mathbb{R}\oplus \mathfrak{h}_{2m+1}$ 
admits a family of complex structures 
$J=J_{(\delta;\, \varepsilon_1, \varepsilon_2,..., \varepsilon_m)}, \delta=c+\im d\in \bC\;( c \not=0), \varepsilon_i=\pm 1$, 
defined by
$$J (T-d Z)=c Z,\; J (c Z)=-(T-d Z),\; J X_i= \varepsilon_i Y_i,\; J Y_i= -\varepsilon_i X_i.$$
Conversely, the above family of complex structures exhausts all
complex structures on $\gg\gh$.
\ep

{\Rm %(i) 
 \label{first:rem} 
In Proposition \ref{heis:prop}, the complex structure $J$ with $J (T-dZ)=cZ,\; J X_i=\varepsilon_i Y$ and
that with $J (T-dZ)=-cZ,\; J X_i=-\varepsilon_iY_i$ are conjugate, and biholomorphic (with the opposite orientations if 
$m$ is even). 
%(ii)  Since $c \not=0$, we can set $c=1$ by taking  a multiple of $T$ as a generator of the center;
%and thus we can set
%\be \label{CT} J T= (-d) T + (1+d^2) Z,\, J Z=-T + d Z,\, J  X_i = \pm Y_i,\,J  Y_i = \mp X_i. \ee
%$$J T= (-d) T + (1+d^2) Z,\, J Z=-T + d Z,\, J  X_i = \varepsilon_iY_i,\,J  Y_i = -\varepsilon_iX_i.$$
}

%%%%%%%%%%%%%%%%%%%%%%%%%  NEWSECTION %%%%%%%%%%%%%%%%%%%%%%%%%%%
%%%%%%%%%%%%%%%%%%%%%%%%%  NEWSECTION %%%%%%%%%%%%%%%%%%%%%%%%%%%

\section{Vaisman structures on $\mathfrak{gl}(2,\bR), \gu(2)$ and $\gg\gh_{2m+2}$} \label{Vaisman_str}

\bd
An {\em almost Hermitian structure}  on a Lie algebra $\gg$ is a pair $(\langle \cdot , \cdot \rangle, J)$ consisting of a 
scalar product $\langle \cdot , \cdot \rangle$ and  a skew-symmetric complex structure $J\in \mathrm{\so}(\gg)$. The triple
$(\gg , \langle \cdot , \cdot \rangle, J)$ is called an {\em almost Hermitian Lie algebra}. The almost Hermitian structure 
is called {\em integrable} (in which case the adjective \emph{almost} is dropped) if the Nijenhuis tensor $N_J$ (\ref{Nj:eq})
vanishes, that is,
%%$$ N_J (X,Y) := [JX,JY] -[X,Y] -J[X,JY]-J[JX,Y] =0$$
$N_J(X,Y)=0$ for all $X,Y\in \gg$.
\ed 

\bd
A Hermitian Lie algebra $(\gg, \langle \cdot , \cdot \rangle, J)$ is called {\em K\"ahler} if its 
fundamental $2$-form $\omega =  \langle \cdot , J \cdot \rangle$ is closed. It is called  {\em locally 
conformally K\"ahler}  (shortly {\em lcK)} if 
$$ d\o =  \omega \wedge \theta$$ 
for some closed $1$-form $\theta \in \gg^*$ (called the {\em Lee form}). An lcK Lie algebra is called {\em Vaisman} 
if $\nabla \theta =0$ holds, where $\nabla \in \gg^*\ot \so (\gg)$ denotes the Levi-Civita connection. Allowing 
possibly indefinite scalar products, we also consider the notion of {\em locally conformally pseudo-K\"ahler}
Lie algebras. Again they are called {\em Vaisman} 
if $\nabla \theta =0$ holds. 
\ed 

Note that lcK (respectively,  Vaisman) Lie algebras $(\gg, \langle \cdot , \cdot \rangle, J)$ correspond to Lie groups $G$ 
with left-invariant lcK (respectively,  Vaisman) structure. The correspondence is one-to-one if we restrict to simply 
connected Lie groups. 

\bd
A {\em contact metric structure} on a Lie algebra $\gg$ of dimension $2n+1$ is a quadruple $(\phi, \eta, \langle \cdot , \cdot \rangle, \widetilde{J})$
consisting of a contact structure $\phi \in \gg^*$, $\phi \wedge (d \phi)^n \not= 0$,\;
$\eta \in \gg$ (called the {\em Reeb field}), $i(\eta) \phi = 1, i(\eta) d \phi = 0$,
a $(1, 1)$-tensor $\widetilde{J} \, , \widetilde{J}^2 = -I + \phi \otimes \eta$
and a scalar product $\langle \cdot , \cdot \rangle$, $\langle X, Y \rangle = \phi(X) \phi(Y) + d \, \phi (\widetilde{J} X, Y)$
for all $X, Y \in \gg$.
\ed

\bd
A {\em Sasaki structure} on a Lie algebra $\gg$ is
a contact metric structure
$(\phi, \eta, \langle \cdot , \cdot \rangle, \widetilde{J})$ satisfying $\langle [\eta, X], Y\rangle + \langle X, [\eta, Y]\rangle= 0$ for
all $X, Y \in \gg$ (Killing field),
and the integrability of  $J = \widetilde{J}|{\cal D}$ on ${\cal D} = {\rm ker} \, \phi$ (CR-structure).
\ed

Note that for any simply connected Sasaki Lie group $G$, its {\em K\"ahler cone} $C(G)$ is defined as
$C(G)=\bR_+ \times G$ with the K\"ahler form $\Omega=r d r \wedge \phi + \frac{r^2}{2}d \phi$,
where a compatible complex structure $\widehat{J}$ is defined by 
$\widehat{J} \eta = \frac{1}{r} \partial_r,\; \widehat{J} \partial_r = (-r) \, \eta$ and $\widehat{J} |{\cal D} = J$.
For any Sasaki Lie group $G$ with contact form $\phi$,  we can define
an lcK form $\omega=\frac{2}{r^2} \Omega = \frac{2}{r} d r \wedge \phi + d \phi$;
or taking $t= -2 \,{\rm log}\, r$,
$\omega= - d t \wedge \phi + d \phi$ on $\bR \times G$,
which is of Vaisman type. We can define a family of complex structures $\widehat{J}$ compatible with
$\omega$ by
$$\widehat{J} (\partial_t-d \eta)=c \eta,\; \widehat{J} (c \eta)=-(\partial_t-d \eta),\; \widehat{J} |{\cal D} = J,$$
where $c, d \in \bR \,(c \not=0))$.
In other words,  we can express $\widehat{J}$ as
$$\widehat{J} \partial_t = \left(-\frac{d}{c}\right)\partial_t + \left(\frac{c^2+d^2}{c}\right) \eta,\; \widehat{J} \eta = 
\left(-\frac{1}{c}\right) \partial_t + \left(\frac{d}{c}\right) \eta,\; \widehat{J} |{\cal D} = J.$$

%$$J \,\partial_t = -d \,\partial_t + (1+d^2) \eta, J \eta = - \partial_t +d \eta,$$
%where $d \in \bR$ and the Lee field is $J \eta$. 

Conversely, any simply connected Vaisman Lie group
is of the form $\bR \times G$ with lcK structure as above,
where $G$ is a simply connected Sasaki Lie group.

\vspace{10pt}

We already know that the only reductive Lie algebras which admit lcK structures are
$\bR \oplus {\gs\gl}(2, \bR)$ and $\bR \oplus {\gs\gu}(2)$ (\cite{ACHK}, \cite{HK1}). 
%Consider lcK structures on the reductive Lie algebras. We know that there are two types:
%$\gg_1=\bR \oplus {\mathfrak sl}(2, \bR)$ and $\gg_2=\bR \oplus {\gs\gu}(2)$. 
For  $\gg_1 = \bR \oplus {\gs\gl}(2, \bR)$, 
take a basis $\{X_1, Y_1, Z_1 \}$ for ${\gs\gl}(2, \bR)$
with Lie brackets defined as in (\ref{sl2:eq})
and $T$ as a generator of the center $\bR$ of $\gg_1$.
Let $t, x_1, y_1, z_1,$ be the Maurer-Cartan forms corresponding to $T, X_1,Y_1,Z_1$
respectively; then we have
$$d t=0, d x_1=y_1 \wedge z_1, d y_1=z_1 \wedge x_1, d z_1=y_1 \wedge x_1$$
and an lcK structure $\omega_1=z_1 \wedge t + y_1 \wedge x_1$
compatible with an integrable complex
structure $J_1$ on $\gg_1$ defined by
$$J_1 X_1=Y_1, J_1 Y_1=-X_1, J_1 T=Z_1, J_1 Z_1=-T.$$

For  $\gg_2 = \bR \oplus \mathfrak{su}(2)$, take a basis $\{X_2, Y_2, Z_2\}$ for $\mathfrak{su}(2)$
with Lie brackets defined as in (\ref{su2:eq})
and $T$ as a generator of the center $\bR$ of $\gg_2$. Let $t, x_2, y_2, z_2,$ be 
the Maurer-Cartan forms corresponding to $T, X_2,Y_2,Z_2$
respectively; then we have
$$d t=0, d x_2=y_2 \wedge z_2, d y_2=z_2 \wedge x_2, d z_2=x_2 \wedge y_2$$
and an lcK structure $\omega_2=z_2 \wedge t + x_2 \wedge y_2$
compatible with an integrable complex
structure $J_2$ on $\gg_2$ defined by
$$J_2 Y_2=X_2, J_2 X_2=-Y_2, J_2 T=Z_2, J_2 Z_2=-T.$$ 

%%%%%%%%%%%%%%%%%%%%%
We know (\cite{ACHK}, \cite{HK2}) that any lcK structure is of the form
$$\omega_{i, \psi} = \psi \wedge t+ d \psi,$$
which is compatible with the above complex structure $J_i$ on $\gg_i$,
where $\psi= a x + b y + c z$ with $a,b,c \in \bR$.

We see that the bilinear form $\langle U,V \rangle _{i, \psi} \,= \omega_{i, \psi}(J_i U,V)$ is
represented, w.r.t. the basis $\{T, X_i, Y_i, Z_i\}$, by the  matrix 
 
$$ A_1= \left(
\begin{array}{cccc}
c & b & -a & 0\\
b & c & 0 & a\\
-a & 0 & c & b\\
0 & a & b & c
\end{array}
\right),\;
A_2= \left(
\begin{array}{cccc}
c & b & -a & 0\\
-b & c & 0 & -a\\
a & 0 & c & -b\\
0 & a & b & c
\end{array}
\right).
$$
Since  $\langle \cdot, \cdot \rangle_{i, \psi}$ must be symmetric and positive definite, we must have
$a=b=0$ for $\gg_2$. 

For $\gg_1$, the characteristic polynomial of $A_1$ is given by
$$\Phi_{A_1} (u) = \{(u-c)^2 -(a^2+b^2)\}^2,$$
and has only positive roots if and only if $c > 0, c^2 > a^2+b^2$. 
The Lee form is $\theta=t$ and the Lee field is
$$\xi=\frac{1}{D} (c T - b X + a Y),$$ 
with $D=c^2 -a^2-b^2$. We have also
$$\langle \xi,\xi\rangle_{1,\psi} \,=\frac{c}{D}.$$
We can see that $\langle [\xi,U],V\rangle _{1,\psi} \,+ \langle U,[\xi,V]\rangle _{1,\psi} \, \not\equiv 0$
unless $a=b=0$. In fact for $U=V=Z_1$, 
$$\langle [\xi, Z_1], Z_1 \rangle _{1,\psi}+\langle Z_1,[\xi,Z_1]\rangle _{1,\psi} \,=2 \langle [\xi, Z_1],
Z_1 \rangle _{1,\psi} \,= -\frac{2}{D}(a^2+b^2)=0$$
if and only if $a=b=0$. 
Conversely for the case $a=b=0$, it is easy to check that 
$$\langle [\xi,U],V\rangle_{1,\psi} \,+ \langle U,[\xi,V]\rangle_{1,\psi} \equiv 0$$
for all $U,V$.
Therefore we have shown the following.

%%%%%%% Proposition %%%%%%%%
\bp 

\begin{enumerate}
\item[(i)] For $J_1$ and $\omega_{1,\psi}$ defined above, $\langle \cdot , \cdot \rangle_{1,\psi}$ defines a 
(positive definite) lcK metric on $\mathfrak{gl}(2,\bR)$ if and only if  $c > 0, c^2 > a^2+b^2$.
It is of Vaisman type if and only if $ c >0,\, a=b=0$;
and it is of non-Vaisman type if and only if $ c>0,\, c^2 > a^2 + b^2 > 0$.

\item[(ii)] For  $J_2$ and $\omega_{2,\psi}$ defined above, $\langle \cdot , \cdot \rangle_{2,\psi}$ defines 
a (positive definite) lcK metric on $\gu(2)$ if and only if $c > 0, a=b=0$; and thus it is always of Vaisman type.
\end{enumerate}
\ep

{\Rm We can see that the Vaisman structure $\omega_1=z_1 \wedge t + y_1 \wedge x_1$
is compatible with the complex structures:
$J_1 X_1=Y_1, J_1 (T-d Z_1)= c Z_1;$
and the Vaisman structure $\omega_2=z_2 \wedge t + x_2 \wedge y_2$ is compatible with the complex structures:
$J_2 Y_2=X_2, J_2 (T- d Z_2)=c Z_2.$
Hence, all complex structures on $\mathfrak{gl}(2,\bR), \gu(2)$
admit their compatible Vaisman structures, up to conjugation (see Remark  \ref{first:rem}).}

\vspace{10pt}

For the case of $\gg\gh_{2m+2}$, which is a nilpotent Lie algebra with a standard basis $\beta =\{T, X_i,Y_i,Z\}$
for which non-zero Lie brackets are given by
$$[X_i,Y_i]=Z,\; i=1, 2, ..., m,$$
We have the canonical Vaisman structure $\omega =z  \wedge t \,+\, \sum_{i=1}^{m} \, y_i \wedge x_i$, which is  compatible with
the complex structure $J=J_{(\delta;\, \varepsilon_1, \varepsilon_2,..., \varepsilon_m)}, \delta=c+\im d\in \bC\;( c \not=0), \varepsilon_i=\pm 1$, 
defined by
$$J (T-d Z)=c Z,\; J (c Z)=-(T-d Z),\; J X_i= \varepsilon_i Y_i,\; J Y_i= -\varepsilon_i X_i,$$
if and only if $\varepsilon_i =1, i=1,2,..., m$.

As in the previous cases, any lcK structure is of the form
$$\omega_\psi=\psi \wedge t + d \psi,$$
where $\psi= \sum_{i=1}^m a_i x_i + \sum_{j=1}^m b_j y_j + c_0 z$
with $a_i, b_j, c_0 \in \bR$. By the compatibility with one of the above complex structures $J$, 
$\langle \cdot, \cdot \rangle_\psi =\omega(J \cdot,\cdot)$ must be symmetric and positive 
definite. Hence $a_i=b_i=0,\; i=1,...,m$, $cc_0>0$ and $c_0\varepsilon_i>0$. 
%As in the previous cases,
%the Vaisman structure $\omega$ is compatible with the complex structures: $J X_i=Y_i, J (T-d Z)= Z,\; i=1,2,...,m.$
%Therefore, we have shown the following.

\bp 
For  $J$ and $\omega_\psi$ defined above, $\langle \cdot , \cdot \rangle_\psi$ defines a (positive definite) lcK metric
on $\gg\gh_{2m+2}$ if and only if $cc_0 > 0$ and $a_i=b_i=0,\;\mathrm{sgn}(c_0)\varepsilon_i=1$ for $i=1,2,..,m$;
and it is always of Vaisman type. More generally, $\langle \cdot , \cdot \rangle_\psi$ is a locally conformally pseudo-K\"ahler 
metric if and only if  $c_0\neq 0$ and $a_i=b_i=0$ for $i=1,2,..,m$; again it is of Vaisman type. 
\ep
%\begin{proof}
{\em Proof.} It remains to check the Vaisman property. A straightforward calculation
shows that the Lee field, the metric dual of the Lee form $\psi$, is of the form $\alpha T +\beta Z$ with 
$\alpha d + \beta =0$ and that this implies $\nabla \psi =0$. 
%(The coefficients of the Lee field can be computed by the second equation $\alpha (cc_0 + \frac{d^2c_0}{c})+\beta \frac{c_0d}{c} =1$ but this is not needed.) 
%\end{proof}
%%%%%%%%%%%%%%%%%%%%%%%%%%%%%  NEWSECTION %%%%%%%%%%%%%%%%%%%%%%%%%%
%%%%%%%%%%%%%%%%%%%%%%%%%%%%%  NEWSECTION %%%%%%%%%%%%%%%%%%%%%%%%%%

\section{Modifications of Hermitian Lie algebras} \label{modification}
In this section we review and discuss in some detail {\em modification} in the category of Hermitian 
and Sasaki Lie algebras.

\vspace{10pt}

%\bd An {\em almost Hermitian structure}  on a Lie algebra $\gg$ is a pair $(\langle \cdot , \cdot \rangle, J)$ consisting of a 
%scalar product $\langle \cdot , \cdot \rangle$ and  a skew-symmetric complex structure $J\in \mathrm{\so}(\gg)$. The triple
%$(\gg , \langle \cdot , \cdot \rangle, J)$ is called an {\cmssl almost Hermitian Lie algebra}. The almost Hermitian structure 
%is called {\cmssl integrable} (in which case the adjective \emph{almost} is dropped) if the Nijenhuis tensor $N_J$ of $J$ 
%vanishes, that is,
%$$ N_J (X,Y) := [JX,JY] -[X,Y] -J[X,JY]-J[JX,Y] =0$$
%for all $X,Y\in \gg$.
%\ed 

Given a Lie group $G$ with $\mathrm{Lie}\, G = \gg$, there is a one-to-one correspondence between 
(almost) Hermitian structures on $\gg$  and left-invariant (almost) Hermitian structures on $G$.  We denote 
by $\so (\gg)$ (respectively, $\gu (\gg)$) the orthogonal (respectively, unitary) Lie algebra of $(\gg , \langle \cdot , \cdot \rangle, J)$;
and by $\mathrm{Der}_{\gu}(\gg)$ the subalgebra of $\gg\gl(\gg)$ consisting of skew-Hermitian derivations of $\gg$: 
$$\mathrm{Der}_{\gu}(\gg) := \mathrm{Der}(\gg)\cap \gu (\gg),$$
where $\mathrm{Der}\,(\gg)$ is the derivation algebra of $\gg$.
 
\bd 
A {\em modification} of an almost Hermitian Lie algebra $(\gg , \langle \cdot , \cdot \rangle , J)$ is a Lie algebra homomorphism 
$$ \phi : \gg \ra \mathrm{Der}_{\gu}(\gg)$$
such that
\begin{enumerate}
\item[(i)] $\phi ([\gg ,\gg ]) =0$ and 
\item[(ii)] $\phi (\mathrm{Im}\, \phi (\gg))=0$, where $\mathrm{Im}\, \phi (\gg) := \phi (\gg ) \gg = 
\{ AX \mid A\in \phi (\gg ), X\in \gg \}$  
denotes the image of the linear Lie subalgebra $\phi (\gg )\subset \mathrm{Der}_{\gu}(\gg)$.
\end{enumerate}
\ed 

%%%%%%% Proposition %%%%%%%%%
\bp 
Given a modification $\phi$ of an almost Hermitian Lie algebra $(\gg , \langle \cdot , \cdot \rangle, J)$,
we define Lie bracket $[\cdot , \cdot ]_\phi$ on $\gg$ by 
$$ [X,Y]_\phi := [X,Y] + \phi (X)Y-\phi(Y)X.$$
Denoting the resulting Lie algebra by $\gg_\phi$, we obtain a new  almost Hermitian Lie algebra 
$(\gg_\phi , \langle \cdot , \cdot \rangle, J)$. Then $(\gg_\phi , \langle \cdot , \cdot \rangle, J)$ is 
integrable if and only if $(\gg , \langle \cdot , \cdot \rangle, J)$ is. 
\ep 
 
%\bpf
{\em Proof.}
Using $\phi ([\gg ,\gg ])=\phi (\mathrm{Im}\, \phi (\gg))=0$ and the abbreviation $\phi_X = \phi(X)$, we compute 
\begin{eqnarray*} [[X, Y]_\phi, Z]_\phi & =& [[X, Y] + \phi_XY - \phi_YX, Z]_\phi\\
&=& [[X, Y] + \phi_XY - \phi_YX, Z] 
- \phi_Z ([X, Y] + \phi_XY - \phi_YX)\\
&=&  [[X, Y] + \phi_XY - \phi_YX, Z] - [\phi_ZX, Y] -[X,\phi_ZY]-\phi_Z\phi_XY +\phi_Z\phi_YX.
\end{eqnarray*}
So 
\begin{eqnarray*}  \sum_{cyclic}[[X, Y]_\phi, Z]_\phi &=& \sum_{cyclic} (-\phi_Z\phi_XY +\phi_Z\phi_YX)=
\sum_{cyclic} (-\phi_X\phi_YZ + \phi_Y\phi_XZ)\\
&=&  
-\sum_{cyclic}  \phi_{[X,Y]}Z =0.
\end{eqnarray*}
Using that $[\phi (\gg ),J]=0$, one can easily check that the Nijenhuis tensor $N_J^\phi$ of $J$ with respect 
to the modified Lie bracket $[\cdot ,\cdot ]_\phi$ coincides with $N_J$ . 
%\epf

%%%%%%% Proposition %%%%%%%%
\bp
\label{modif_lcK:prop} Let  $(\gg_\phi , \langle \cdot , \cdot \rangle, J)$ be a modification of a Hermitian
Lie algebra $(\gg , \langle \cdot , \cdot \rangle, J)$. Then $(\gg_\phi , \langle \cdot , \cdot \rangle, J)$ K\"ahler (respectively, lcK) 
if and only if $(\gg , \langle \cdot , \cdot \rangle, J)$ is K\"ahler (respectively, lcK). 
\ep

{\em Proof.}
The Levi-Civita connection $\n^\phi$ of $(\gg_\phi, \langle \cdot , \cdot \rangle)$ is given related to the Levi-Civita connection
$\n$ of $\gg$ by 
\be \label{LC:eq} \n^\phi = \n + \phi.\ee
In fact,  first $\n^\phi$ is metric since $\n$ is metric and $\phi\in \gg^*\ot \so(\gg)$. Second, the torsion $T^\phi$ of $\n^\phi$ 
is related to the torsion $T$ of $\n$ by 
\begin{eqnarray*}  T^\phi(X,Y) &=& \n^\phi_XY-\n^\phi_YX -[X,Y]_\phi = T(X,Y) +(\phi_XY-\phi_YX) -(\phi_XY-\phi_YX)\\
&=&T(X,Y)=0.
\end{eqnarray*}
Next, since $\phi$ takes values in the unitary Lie algebra $\gu (\gg)=\{ A\in \so (\gg) \mid [A,J] =0\}$
it follows from \re{LC:eq} that 
$$ \n^\phi J = \n J\quad\mbox{and}\quad \n^\phi \o = \n \o.$$ 
As a consequence, $d^\phi \omega =  d\omega$, where $d^\phi$ denotes the (Chevalley-Eilenberg)
differential  in the Lie algebra $\gg_\phi$. 
This proves that $\o$ is closed (respectively, conformally closed) in $\gg_\phi$ if and only if it is in $\gg$. 
%\epf

%%%%%%% Proposition %%%%%%%%
\bp
\label{VaisProp} A  modification $(\gg_\phi , \langle \cdot , \cdot \rangle, J)$ of a Vaisman Lie algebra $(\gg , \langle \cdot , \cdot \rangle, J)$ 
(of dimension $\ge 4$) with Lee form $\theta$ is Vaisman 
if and only if $\theta (\mathrm{Im}\,\phi (\gg))=0$. 
\ep 

{\em Proof.}
By Proposition \ref{modif_lcK:prop} we know that $(\gg_\phi , \langle \cdot , \cdot \rangle, J)$ is lcK. It is Vaisman if and only if 
$\n^\phi \theta =0$. Since $\n \theta =0$, the latter condition amounts to $\theta \circ \phi_X =0$ for all $X\in \gg$, which is 
equivalent to $\theta (\mathrm{Im}\,\phi (\gg))=0$.  
%\epf

%%%%%%%%%%%%%%%%%%%%%%%%%%%%%%  NEWSECTION %%%%%%%%%%%%%%%%%%%%%%%%%%
%%%%%%%%%%%%%%%%%%%%%%%%%%%%%%  NEWSECTION %%%%%%%%%%%%%%%%%%%%%%%%%%

\section{Unimodular Sasaki and Vaisman Lie algebras} \label{unimodular}
  
Let us denote by  $\gt^q \cong \mathbb{R}^q$ the Cartan subalgebra $\gt^q \subset \gu(q)$ consisting of purely imaginary diagonal matrices. 
Given a linear map
$$ \psi : \mathbb{R}^{2p+1} \ra \gt^q,$$
we can define a meta-abelian Lie algebra 
$\gm :=  \mathbb{R}^{2p+1} \ltimes_\psi \mathbb{C}^q$ by 
considering  $\mathbb{R}^{2p+1}$ as an abelian subalgebra acting on the 
abelian ideal $\mathbb{C}^q$ by $\psi$. We decompose $\mathbb{R}^{2p+1} = \mathbb{R}Z_0 \oplus \mathbb{C}^p$ and 
consider the standard symplectic form $\o$ on $V= \mathbb{C}^p \oplus \mathbb{C}^q = \mathbb{C}^m= \mathbb{R}^{2m}$, $m=p+q$. 

\bl 
The pull back of $\o$ via the projection $\gm  \ra V, \lambda Z_0 + v\mapsto v$ ($\lambda\in \mathbb{R}, v\in V$), 
defines a $2$-cocycle $\o_\gm$ in $\gm$. 
\el 
 
{\em Proof.}
We check that $\sum_{cyclic}\o_\gm ([X,Y],Z)=0$ for all $X,Y,Z\in \gm$. This is obviously satisfied if 
all three vectors belong to the abelian subalgebra $\mathbb{R}Z_0 \oplus \mathbb{C}^p$ or to the 
abelian ideal $\mathbb{C}^q$. If $X_1, X_2 \in \mathbb{R}Z_0 \oplus \mathbb{C}^p$ and $Y\in \mathbb{C}^q$
then 
$$ \sum_{cyclic}\o_\gm ([X_1,X_2],Y) =\o_\gm (\psi(X_2)Y,X_1) - \o_\gm (\psi (X_1)Y,X_2)=0,$$
since $\mathbb{C}^p$ and $\mathbb{C}^q$ are $\omega$-orthogonal. Similarly, taking $X\in \mathbb{R}Z_0 \oplus \mathbb{C}^p$ and 
$Y_1,Y_2 \in \mathbb{C}^q$ we obtain 
$$ \sum_{cyclic}\o_\gm ([X,Y_1],Y_2) =  \o_\gm (\psi(X)Y_1,Y_2) -\o_\gm (\psi (X)Y_2,Y_1)=0,$$
since $\psi(X)\in \gu(q) \subset \gsp (\mathbb{R}^{2q})$. 
%\epf
 
\noindent 
Next we  consider  the one-dimensional central extension 
$$ 0\ra \mathbb{R}Z_1 \ra \gg\gh (\psi) \ra \gm \ra 0$$ 
of $\gm$ by the 
cocycle $\o_\gm$.  Explicitly, this means that the Lie bracket of $\gg\gh (\psi)$ is related to the Lie bracket $[\cdot ,\cdot ]_\gm$ of $\gm$  by 
$$ [X,Y] = [X,X]_\gm +\o_\gm (X,Y)Z_1,\quad X,Y\in \gm,$$ 
and $[Z_1,\gm]=0$.  Note that $\gg\gh (\psi)$ is a solvable Lie algebra, which is at most of $3$-step type. 
 
\bt 
\label{mainThm} Let $\gg$ be a unimodular  Lie algebra (of dimension $\ge 4$) which admits a Vaisman structure. 
Then $\gg$ is isomorphic to one of the following. 
\begin{enumerate}
\item[(i)] $\gg = \gu(2) = \mathbb{R} \oplus \su (2)$, 
\item[(ii)] $\gg = \ggl(2,\bR)=\mathbb{R} \oplus \gsl (2,\mathbb{R})$, or 
\item[(iii)] a solvable unimodular Lie algebra $\gg\gh (\psi)$ associated with a linear map $\psi : \mathbb{R}^{2p+1} \ra \mathbb{R}^q$ as 
described above. 
\end{enumerate}
%(See Remark \ref{rem1} below for a description of all Vaisman structures on these Lie algebras)
\et 
  
{\em Proof.}
By \cite[Theorem 2.1]{AHK} $\gg$  is a modification of one of the following Lie algebras endowed with a Vaisman structure: $\gu (2)$, 
$\ggl (2,\mathbb{R})$, or $\gg\gh=\mathbb{R} \oplus \mathfrak{h}_{2m+1}$.
%where $\mathfrak{h}_{2m+1}$ denotes the Heisenberg Lie algebra of dimension $2m+1\ge 3$
We claim that in the first two cases any modification $\gg_\phi$ of $\gg$ is isomorphic to $\gg$. 
This follows from the next lemma. To state it we first need to define the notion of a modification for Riemannian (rather than almost Hermitian) 
Lie algebras. 

\bd
A Riemannian Lie algebra  $(\gg , \langle \cdot, \cdot \rangle)$ is a Lie algebra endowed with a scalar product. A {\em modification} of
$(\gg , \langle \cdot, \cdot \rangle)$ is a Lie algebra homomorphism 
$$ \phi : \gg \ra \mathrm{Der}_{\so}(\gg) 
= \mathrm{Der}(\gg)\cap \so (\gg)$$
such that $\phi ([\gg ,\gg ]) =0$ and $\phi (\mathrm{Im}\, \phi (\gg))=0$.
\ed 
As before, a modification $\phi$ of a Riemannian Lie algebra $(\gg , \langle \cdot, \cdot \rangle)$ gives rise to a
new Riemannian Lie algebra $(\gg_\phi, \langle \cdot, \cdot \rangle)$.
The Lie algebra $\gg_\phi$ as well as the Riemannian Lie algebra $(\gg_\phi, \langle \cdot, \cdot \rangle)$ 
will be again called a modification of $(\gg , \langle \cdot, \cdot \rangle)$. 

\bl
Let $\gg = \mathbb{R}\oplus \gs$, $\gs = [\gg , \gg]$, be a reductive Lie algebra with one-dimensional center 
and $\langle \cdot, \cdot \rangle$ a scalar product on $\gg$. Then any modification
$\gg_\phi$ of $(\gg , \langle \cdot, \cdot \rangle)$ is isomorphic to $\gg$. 
\el 

{\em Proof.}
Let $Z$ be a generator of the center of $\gg$. Any modification $\phi : \gg \ra \mathrm{Der}_{\so}(\gg)$ is
completely determined by $\phi (Z)\in \mathrm{Der}_{\so}(\gg)$, since $\phi ([\gg , \gg]) = \phi (\gs )=0$. 
Since the center is invariant under any derivation and is one-dimensional, any skew-symmetric derivation necessarily 
acts trivially on the center. This proves that $\mathrm{Der}_{\so}(\gg)$ consists of inner derivations, that is 
$\phi (Z) = \mathrm{ad}_{X_0}$ for some element $X_0\in \gs$. The modified Lie bracket is therefore 
given by 
$$ [X+ \lambda Z,Y+\mu Z]_\phi = [X,Y] + \lambda\, \mathrm{ad}_{X_0}Y-\mu\, \mathrm{ad}_{X_0}X,
\quad X,Y \in \gs,\quad \lambda, \mu \in \mathbb{R}.$$
One can easily check that the map 
$$ X+\lambda Z \mapsto X+\lambda Z +\lambda X_0\quad (X\in \gs, \lambda \in \mathbb{R})$$
defines an isomorphism $\gg_\phi \cong \gg$. 
%\epf

Next we will describe all possible modifications of the Vaisman Lie algebra 
$ \gg\gh = \mathbb{R} \oplus \mathfrak{h}_{2m+1}.$
We have seen in Section \ref{Vaisman_str} that (up to an automorphism of $\gg\gh$) 
the Vaisman structure $(\langle \cdot , \cdot \rangle, J)$ is given by the canonical Hermitian structure
on $V=\mathbb{C}^m=\mathbb{R}^{2m}$ which is extended
to $\gg\gh = \mathbb{R}Z_0 + \mathbb{R}Z_1 + V$ as follows, where $Z_1$ denotes a generator of 
the center of $\gh_{2m+1}=  \mathbb{R}Z_1 + V$ such that 
$[X,Y] = \o (X,Y)Z_1$ for all $X, Y\in V$ and $Z_0$ denotes a generator of the $\mathbb{R}$-factor in
$\gg\gh = \mathbb{R} \oplus \mathfrak{h}_{2m+1}$. 
The scalar product of $V$ is extended such that $Z_0$, $Z_1$ are orthonormal and perpendicular to $V$.
The complex structure of $V$ is extended such that $JZ_0=Z_1$. The Lee form $\theta$ is then given by
$\theta = \langle Z_0, \cdot \rangle$. We denote this Vaisman structure by 
$(\langle \cdot ,\cdot \rangle, J)$. 

\vspace{10pt}

Using the invariance of the center and the derived ideal under any derivation, it is easy to see that
$\mathrm{Der}_{\gu}(\gg\gh) = \gu (m)$ with the natural action on $V=\bC^m\subset \gg\gh$.
Let $\phi : \gg\gh \ra \gu(m)$ be a modification of the Vaisman Lie algebra 
$(\gg\gh ,  \langle \cdot, \cdot \rangle, J)$. Note that due to 
$\phi (\gg\gh ) \subset \gu (m) \subset \ker \theta$, $(\gg\gh_\phi, \langle \cdot ,\cdot \rangle, J)$ 
is again a Vaisman Lie algebra in virtue of Proposition \ref{VaisProp}.  
Since $\phi ([\gg\gh , \gg\gh ])=0$, $\phi$ factorizes through the abelian quotient Lie algebra 
$\gg\gh/\mathbb{R}Z_1\cong \mathbb{R} Z_0 \oplus V$.  In particular, $\phi (\gg\gh) \subset \gu(m)$ is
abelian and we can assume 
(up to conjugation in $\gu(m)$) that $\phi$ takes values in the standard Cartan subalgebra $\gt^n \subset \gu (m)$. 
We denote by $\bar{\phi} : \mathbb{R} Z_0 \oplus V \ra \gt^m$ the induced homomorphism. Since $\phi (\gg\gh )$ is a subalgebra 
of $\gu (\gg\gh )$, we can decompose $V=V_0\oplus V_1$ as an orthogonal sum of $J$-invariant subspaces, 
where $V_0 = \ker \phi (\gg\gh) := \{ v\in V \mid \phi (\gg\gh ) v =0\}= \bC^p$ and $V_1 = \mathrm{Im}\, \phi (\gg\gh)=\bC^q$.  
This implies that $\phi : \gg\gh \ra \gt^q\subset \gu(V_1)=\gu(q)$.
Next the condition $\phi (\mathrm{Im}\, \phi( \gg\gh))$ implies that 
$\bar\phi$ vanishes on $V_1$ and, hence, can be considered as a linear map $\psi : \mathbb{R} Z_0 \oplus V_0 \ra \gt^q$. 
Now one can check that the modification $\gg\gh_\phi$ coincides with the Lie algebra $\gg\gh(\psi)$ described above. 
Note that $\gg\gh (\psi)$ is a unimodular solvable Lie algebra since modifications preserve unimodularity.
%The unimodularity of 
%$\gg\gh (\psi)$ follows,  for example, from the fact that nilpotent Lie algebras are unimodular and modifications preserve unimodularity. 
%\epf

\bc
Let $\gh$ be a unimodular  Lie algebra which admits a Sasaki structure. Then $\gg$ is isomorphic to one of the 
following. 
\begin{enumerate}
\item[(i)] $\gh = \su (2)$, 
\item[(ii)] $\gh = \gsl (2,\mathbb{R})$, or 
\item[(iii)] a solvable unimodular Lie algebra $\gh_{2m+1} (\varphi)$ associated with a linear map
$\varphi : \mathbb{R}^{2p} \ra \mathbb{R}^q$ as described above. 
\end{enumerate}
%{\color{red} Include a remark specifying all Sasaki structures on the above groups.}
\ec 

 %%%%%%%%%%%%%%%%%%%%%%%%%%%%  NEWSECTION %%%%%%%%%%%%%%%%%%%%%%%%%%
  %%%%%%%%%%%%%%%%%%%%%%%%%%%%  NEWSECTION %%%%%%%%%%%%%%%%%%%%%%%%%%
 
\section{Biholomorphism types of $\mathfrak{gl}(2,\bR), \gu(2)$ and $\gg\gh_{2m+2}$} \label{biholomorphism}

In this section, we determine the biholomorphism type of each simply connected unimodular Lie group
listed in the main theorem; namely we show that $\bR \times \widetilde{SL}(2, \bR)$, $\bR \times {\rm SU}(2)$, 
and $GH_{2m+2}=\bR \times  H_{2m+1}$ with its modification $GH_{2m+2}(\psi)$
are biholomorphic to $\bH \times \bC$, $\bC^2\backslash \{0\}$,  and $\bC^{m+1}$ respectively.
Note that this result is concerned with all invariant complex structures on these Lie groups we have determined in Section \ref{complex_str},
which is independent of the existence of Vaisman structures on them. It turned out as seen in Section \ref{Vaisman_str} 
that all these complex structures admit compatible Vaisman structures.

%%%%%%%%%%%%%%%%%%%  $\bR \times \widetilde{SL}(2, \bR)$
\vspace{10pt}

For the case of $\bR \times \widetilde{SL}(2, \bR)$, we first consider a diffeomorphism $\Phi$:

$$\Phi: GL^+(2, \bR) \longrightarrow \bH \times \bC^*$$
defined by
$$ g = \left(
\begin{array}{cc}
a & b\\
c & d
\end{array} 
\right)
\longrightarrow ( \left( \frac{ac +bd}{c^2+d^2} \right) + \sqrt{-1} \left( \frac{D}{c^2+d^2} \right),\; d + \sqrt{-1}\, c\,), $$
where $D = {\rm det}\, g = ad-bc > 0$. Note that by the Iwasawa decomposition $g$ can be expressed uniquely as
$$g = \sqrt{D}
\left(
\begin{array}{cc}
1 & x\\
0 & 1
\end{array}
\right)
\left(
\begin{array}{cc}
\sqrt{y} & 0\\
0 & \sqrt {y}^{-1}
\end{array} 
\right)
\left(
\begin{array}{cc}
\cos \theta & -\sin \theta\\
\sin \theta & \cos \theta
\end{array} 
\right),
$$
where 
$$x =  \left( \frac{ac +bd}{c^2+d^2} \right),\; y =  \left( \frac{D}{c^2+d^2} \right),\; e^{\sqrt{-1}\, \theta} = \frac{d + \sqrt{-1}\, c}{\sqrt{c^2+d^2}}.$$
If we consider the canonical action of $GL^+(2, \bR)$ on $\bH$, then we have $g \cdot \sqrt{-1} = x+\sqrt{-1}y$ with the isotropy
subgroup $\bR^+ \times SO(2)$ at $\sqrt{-1}$. We see that $\Phi$ defines a biholomorphism between $GL^+(2, \bR)$ and $\bH \times \bC^*$.
In fact, recall that for  $\gg_1 = \bR \oplus \mathfrak{sl}(2, \bR)$, 
we have a basis $\{T, X_1, Y_1, Z_1 \}$ of $\gg_1$ as
$$
T =\frac{1}{2}\left(
\begin{array}{cc}
1 & 0\\
0 & 1
\end{array}
\right),
X_1=\frac{1}{2}\left(
\begin{array}{cc}
0 & 1\\
1 & 0
\end{array}
\right),
\;
Y_1=\frac{1}{2}\left(
\begin{array}{cc}
1 & 0\\
0 & -1
\end{array}
\right),
\;
Z_1=\frac{1}{2}\left(
\begin{array}{cc}
0 & -1\\
1 & 0
\end{array}
\right).
$$
They are  pushed forward by $\Phi$
to the vector fields, which are expressed in the local coordinates
$(x + \sqrt{-1} \, y, r e^{\sqrt{-1}\, \theta}),\;
y>0, r = \sqrt{c^2+d^2} >0$ of $\bH \times \bC^*$, as
$$
X_1' = y \cos 2 \theta \frac{\partial}{\partial x}  - y \sin2 \theta \frac{\partial}{\partial y} + 
\frac{1}{2} \cos 2\theta \frac{\partial}{\partial \theta} +\frac{r}{2} \sin 2\theta \frac{\partial}{\partial r},\;
$$
$$
Y_1' = y \sin 2 \theta \frac{\partial}{\partial x} +  y \cos 2 \theta \frac{\partial}{\partial y} + 
\frac{1}{2} \sin 2\theta \frac{\partial}{\partial \theta} - \frac{r}{2} \cos 2 \theta \frac{\partial}{\partial r},\;
$$
$$
Z_1' =  \frac{1}{2} \frac{\partial}{\partial \theta},\;
T' = \frac{r}{2} \frac{\partial}{\partial r},
$$
respectively. For the canonical complex structure $J$ on $\bH \times \bC^*$ defined by
$$J \frac{\partial}{\partial x} = \frac{\partial}{\partial y},\; J \,(r \frac{\partial}{\partial r})= \frac{\partial}{\partial \theta},$$
we have the compatibility 
$$J X_1' = Y_1', J T'=Z_1';$$
and thus $\Phi$ defines a biholomorphic map. Then it induces a biholomorphic map $\bar{\Phi}$
between their universal coverings $(\bR \times \widetilde{SL}(2, \bR), \bar{J})$
and $\bH \times \bC$, where $\bar{J}$ is the induced integrable complex structure from $J$.

\vspace{10pt}

We now consider a biholomorphic map $\Phi_\delta$ from 
$(\bR \times {SL}(2, \bR), J)$ to $\bH \times \bC^*$,
where $J (T-l Z_1) = k Z_1, J X_1 = Y_1$ for $\delta=k + \sqrt{-1} l, k \not=0$.

We consider $$\Phi_\delta: \bR \times {SL}(2, \bR) \longrightarrow \bH \times \bC^*$$
defined by
$$ g = (t, \left(
\begin{array}{cc}
a & b\\
c & d
\end{array} 
\right))
\longrightarrow ( \left( \frac{ac +bd}{c^2+d^2} \right) + \sqrt{-1} \left( \frac{D}{c^2+d^2} \right),\; e^{\delta\, t} \,(d + \sqrt{-1}\, c) ), $$
where $\delta=k+\sqrt{-1}\, l$. Then $X_1'$ and $Y_1'$ are the same as before, whereas $Z_1'$ and $T'$ are changed to the following.
$$
Z_1' =  \frac{1}{2} \frac{\partial}{\partial \theta},\;
T' = k\, \frac{r}{2} \frac{\partial}{\partial r} +  l\,\frac{1}{2} \frac{\partial}{\partial \theta}.
$$

Clearly we have
$$J X_1' = Y_1',\; J (T'-l Z_1') =  kZ_1',$$

and thus $\Phi_\delta$ defines a biholomorphic map.

\vspace{10pt}

%%%%%%%%%%%%%%%%%%%%%%%%%%  $\bR \times SU(2)$

For the case of $\bR \times SU(2)$, 
we consider a canonical diffeomorphism $\Phi_\delta$:

$$\Phi_\delta: \bR \times SU(2) \longrightarrow \bC^2 \backslash \{0\}$$
defined 
$$(t,z_1,z_2) \longrightarrow (e^{\delta \,t} z_1, e^{\delta \,t} z_2),$$

\noindent where $\delta=k+{\im l}$ and
$SU(2)$ is identified with
$$S^3=\{(z_1,z_2) \in \bC^2 \,|\;|z_1|^2+|z_2|^2=1\}$$
by the correspondence
$$
\left(
\begin{array}{cc}
z_1 & -\overline{z}_2\\
z_2 & \overline{z}_1
\end{array}
\right)
\longleftrightarrow (z_1,z_2).
$$
Then $\Phi_\delta$ induces a biholomorphism
between $\bR \times SU(2)$ with $J_\delta$ and 
$\bC^2 \backslash \{0\}$.
In fact, recall that for  $\gg_2 := \bR \oplus \mathfrak{su}(2)$, we set  a basis $\{X_2, Y_2, Z_2 \}$ of  $\mathfrak{su}(2)$ as
$$ 
X_2=\frac{1}{2}\left(
\begin{array}{cc}
0 & -1\\
1 & 0
\end{array}
\right),
\;
Y_2=\frac{1}{2}\left(
\begin{array}{cc}
0 & \sqrt{-1}\\
\sqrt{-1} & 0
\end{array}
\right),
\;
Z_2=\frac{1}{2}\left(
\begin{array}{cc}
\sqrt{-1} & 0\\
0 & -\sqrt{-1}
\end{array}
\right),
$$
and $T$ as a generator of the center $\bR$ . They are  pushed forward by $\Phi_\delta$
to the vector fields given as
$$
X_2' = \frac{1}{2} \,(\,-\bar{z_2} \frac{\partial}{\partial z_1} + \bar{z_1} \frac{\partial}{\partial z_2}),\;
Y_2'= \frac{\sqrt{-1}}{2} \,(\,-\bar{z_2} \frac{\partial}{\partial z_1} + \bar{z_1} \frac{\partial}{\partial z_2}),
$$
$$
Z_2'= \frac{\sqrt{-1}}{2} \,( \,z_1 \frac{\partial}{\partial z_1}  + \,z_2 \frac{\partial}{\partial z_2} -
\, \bar{z}_1 \frac{\partial}{\partial \bar{z}_1}  - \, \bar{z}_2 \frac{\partial}{\partial \bar{z}_2}),\;
T' = \frac{\delta}{2} \,( \,z_1 \frac{\partial}{\partial z_1}  + \,z_2 \frac{\partial}{\partial z_2}) +
\frac{\bar{\delta}}{2} \,( \, \bar{z}_1 \frac{\partial}{\partial \bar{z}_1}  + \, \bar{z}_2 \frac{\partial}{\partial \bar{z}_2}),
$$
respectively. Then we have
$$J X_2' = Y_2',\; J (T'-l Z_2') = k Z_2',$$
and thus $\Phi_\delta$ defines a biholomorphic map.
%\noindent $\Phi_\delta$ induces a biholomorphism
%between $G=S^1 \times SU(2)$ with $J_\delta$ and a primary Hopf surface 
%$S_{\lambda_\delta}=W/\Gamma_{\lambda_\delta}$.

\vspace{10pt}

%%%%%%%%%%%%%%%%%%%% $GH_{2m+2}=\bR \times H_{2m+1}$

We now consider the case of $GH_{2m+2}=\bR \times H_{2m+1}$.
Let $\gh_{2m+1}$ be the Heisenberg Lie algebra of dimension $2m+1$ and $H_{2m+1}$
be the corresponding Lie group. $\gg\gh_{2m+2} = \bR \oplus \gh_{2m+1}$
has a canonical basis $\{X_i, Y_j, Z, T\}$ for which non-zero Lie brackets are defined by
$$[X_i, Y_i] = Z, $$
where $i = 1, ... , m$; and the complex structure $J=J_\delta$, $\delta=k+l\sqrt{-1}$,  defined by\
$$J X_i = \varepsilon_i Y_i,\, J Y_i = - \varepsilon_i X_i,\, J (T-l Z) = kZ,\, J (kZ) = - (T-lZ),$$
where $\varepsilon_i=\pm 1, i =1, ...,m$.

We consider the bijective map $\Phi_\delta: GH_{2m+2} = H_{2 m+ 1} \times \bR \rightarrow \bC^{m+1}$ defined by

$$\big( \left(
\begin{array}[c]{ccc}
1 & {\bf x} & z\\
0 & {\rm I}_m & {\bf y}^t\\
0 & 0 & 1
\end{array}
\right),\, t \,\big)
\rightarrow
({\bf x} + \sqrt{-1} {\bf y},\, (2 kt + \frac{1}{2}(\|{\bf x}\|^2 + \|{\bf y}\|^2)) + \sqrt{-1}\,
(2 (lt+z) - {\bf x} \cdot {\bf y})),
$$
where ${\bf x}=(x_1, x_2, ..., x_m), \, {\bf y}=(y_1, y_2, ..., y_m) \in \bR^m$, ${\bf x} \cdot {\bf y} =
\sum_{i=1}^{m} \varepsilon_i \, x_i y_i$ and $\|{\bf x}\|^2 = {\bf x} \cdot {\bf x}$.

The basis $\{X_i, Y_j, Z, T\}$ of $\gg\gh_{2m+2}$ is given in the coordinates of 
$H_{2 m+ 1} \times \bR$ as
$$X_i= \frac{\partial}{\partial x_i},\; Y_j=\frac{\partial}{\partial y_j} + x_j\, \frac{\partial}{\partial z},\;
Z= \frac{\partial}{\partial z},\; T= \frac{\partial}{\partial t};$$
and they are pushed forward by the map $\Phi_\delta$ to
$$X_i'= \frac{\partial}{\partial x_i} + x_i  \frac{\partial}{\partial t} - y_i  \frac{\partial}{\partial z},\;
Y_j'= \frac{\partial}{\partial y_j}  + y_j  \frac{\partial}{\partial t} + x_j  \frac{\partial}{\partial z},\;
Z'=2\, \frac{\partial}{\partial z},\; T'=2k\frac{\partial}{\partial t}+2l \frac{\partial}{\partial z},$$
respectively. Then we have  
$$J X_i' = Y_i',\; J (T'-l Z') = k Z',\, i=1,2,..,m,$$
and thus $\Phi_\delta$ defines a biholomorphic map (cf. \cite{H}).

\vspace{10pt}

We can express the group operation on $GH_{2m+2}$ in the coordinates of $\bC^{m+1}$ as
$$ ({\bf w}, v) \cdot ({\bf z}, u) = ({\bf w} + {\bf z}, v + \bar{\bf w} \cdot {\bf z} + u),$$
where ${\bf w} = (w_1, w_2, \cdots, w_m), {\bf z} = (z_1, z_2, \cdots, z_m) \in \bC^m$
and $v, u \in \bC$.  %Note that for $\lambda_i, i=1,.., m$ with
%$|\lambda_i| = 1$, we have an automorphism $\lambda$ of $GH_{2m+2}$ which is defined in $\bC^{m+1}$ by
%$$\lambda:  (z_1, z_2, \cdots, z_m, u) \rightarrow (\lambda_1  z_1, \lambda_2 z_2, \cdots, \lambda_m  z_m, u).$$

%%%%%%%%%%%%%%%%%%%%%%%%%%%
\vspace{10pt}

We finally consider the case of the modifications $GH_{2m+2}(\psi)$ of $GH_{2m+2}$.
Recall (see Section \ref{modification}) that $\gg\gh(\psi)$ is expressed as the central extension of $\gm$, where
$\gm = \bR^{2p+1} \ltimes_\psi \, \bC^{q}$ with the action $\psi: \bR^{2p+1} \rightarrow \gt$,
where $\gt$ is the cartan subalgebra of $\gu(q)$. Since $\gt$ is abelian subalgebra of $\gu(q)$,
we can take a canonical basis $\beta=\{X_i, Y_j, Z_0\}$ of $\gm$ for which
non-zero Lie brackets is given by
$$[X_i, X_{p+j}] = Y_{p+j},\;  [X_i, Y_{p+j}] = -X_{p+j},$$
$$[Y_i, X_{p+j}] = Y_{p+j}, \; [Y_i, Y_{p+j}] = -X_{p+j},$$
$$[Z_0, X_{p+j}] = Y_{p+j},\;  [Z_0, Y_{p+j}] = -X_{p+j},$$
where $p+q=2m$, $i=1,2,...,p, j=1,2,...,q$. Note that  some of the adjoint actions ${\rm ad}\, X_i, {\rm ad}\, Y_j$ or
${\rm ad}\, Z_0$ can be trivial; the case all of them are trivial corresponds the case $\gm=\bR \times \bC^m$.
 %$ i=1, 2,..., p, j \in \Lambda_i \subset \{1, 2,...,q\}.$
Then $\gg\gh_{2m+2}(\psi)$ is a central extension of $\gm$ with the generator $Z_1$ of $\bR$
$$0 \rightarrow \bR Z_1 \rightarrow \gg\gh_{2m+2}(\psi) \rightarrow \gm \rightarrow 0,$$
where the additional non-zero Lie brackets is given by
$$[X_i, Y_i]=Z_1,\; i=1,2,...,m.$$
We can see that the complex structure $J$ on $\gg\gh_{2m+2}$:
$$J X_i=Y_i, J Y_j=-X_j, J Z_0=Z_1, J Z_1=-Z_0,\; i, j=1,2,...,m.$$
is also integrable for $\gg\gh_{2m+2}(\psi)$.
 
\vspace{10pt}

Let $GH_{2m+2}(\psi)$ be the Lie group corresponding to $\gg\gh_{2m+2}(\psi)$, which is a central extension of 
$\bR^{2p+1} \ltimes \, \bC^{q}$:
$$1 \rightarrow \bR \rightarrow GH_{2m+2}(\psi) \rightarrow  \bR^{2p+1} \ltimes \, \bC^{q} \rightarrow 1.$$
where the group structure on $\bR^{2p+1} \ltimes \, \bC^{q}$ is given by the action
$\psi: \bR^{2p+1} \rightarrow {\rm Aut}(\bC^q)$:
$$\psi(\bar{t}_i) ((z_1, z_2, \ldots, z_q)) = 
(e^{\sqrt{-1}\, a^i_1\, t_i} z_1, e^{\sqrt{-1}\, a^i_2\, t_i} z_2, ...,
e^{\sqrt{-1}\, a^i_q\, t_i} z_q),$$
with $\bar{t}_i= t_i e_i$ ($e_i$: the $i$-the unit vector in $\bR^{2p+1})$, and
$a^i_j \in \bR$, $i=1, 2,..., 2p+1$. To be more precise its group multiplication is given as
$$({\bf t}, w_1, w_2,..., w_q) \cdot ({\bf s}, z_1, z_2,..., z_q) = 
({\bf t}+{\bf s}, w_1 + e^{\sqrt{-1}\, {\bf a}_1 \cdot {\bf t}} z_1,\,  e^{\sqrt{-1}\, {\bf a}_2 \cdot {\bf t}} z_2,...,\,  e^{\sqrt{-1}\, {\bf a}_q \cdot {\bf t}} z_q),$$
where ${\bf t}, {\bf s}, {\bf a}_j \in \bR^{2p+1},\; j=1,2,..., q$. This group multiplication can be extended holomorphically to
$GH_{2m+2}(\psi)$, since the map $\lambda$ defined on $GH_{2m+2}$ in the coordinates $\bC^{m+1}$:
$$\lambda:  (z_1, z_2, \cdots, z_m, u) \rightarrow (\lambda_1  z_1, \lambda_2 z_2, \cdots, \lambda_m  z_m, u),$$
for  $|\lambda_i|=1, i=1,2,...,m$, is an automorphism of $GH_{2m+2}(\psi)$.
%$ ({\bf w}, v) \cdot ({\bf z}, u) =$
%$$(w_1 + e^{2 \pi \eta_1 t}  z_1, w_2 + e^{2 \pi \eta_2 t}  z_2, \cdots,
%w_m + e^{2 \pi \eta_m t}  z_m,\; v +  \bar{\bf w} \cdot {\bf z} + u),$$
%where $\eta_i \in \bR$ and $t = {\rm Re}(v)$.
In particular, $GH_{2m+2}(\psi)$ is biholomorphic to $\bC^{m+1}$.

%%%%%%%%%%%%%%%%%%%%%%%%%%%%%% REFERENCES %%%%%%%%%%%%%%%%%%%%%%
%%%%%%%%%%%%%%%%%%%%%%%%%%%%%% REFERENCES %%%%%%%%%%%%%%%%%%%%%%

\address{Vicente Cort\'es\\
Department Mathematik\\
und Zentrum f\"ur Mathematische Physik\\
Universit\"at Hamburg, Bundesstra{\ss}e 55\\
D-20146 Hamburg, Germany}
{cortes@math.uni-hamburg.de}
\address{Keizo Hasegawa\\
Department of Mathematics\\
Graduate School of Science\\
Osaka University, Machikaneyama-cho\\
ToyonakaOsaka 560-0043, Japan}
{hasegawa@math.sci.osaka-u.ac.jp}
\vspace{-10pt}

\hspace{7.38cm}
\address{
Department of Mathematics\\
Faculty of Education\\
Niigata University\\
8050 Ikarashi-Nino-cho\\
Niigata 950-2181, JAPAN}
{hasegawa@ed.niigata-u.ac.jp}

\end{document}